\documentclass[a4paper, 12pt]{amsart}
\usepackage{amsmath}
\usepackage{amssymb}
\usepackage{amsthm}
\usepackage{amscd}
\usepackage[all]{xy}

\theoremstyle{plain}
\newtheorem{thm}{Theorem}[section]
\newtheorem{prop}[thm]{Proposition}
\newtheorem{cor}[thm]{Corollary}
\newtheorem{lem}[thm]{Lemma}

\newtheorem{mainthm}{Theorem}
\newtheorem{maincor}{Corollary}
\theoremstyle{definition}
\newtheorem{defn}[thm]{Definition}
\newtheorem{expl}[thm]{Example}
\theoremstyle{remark}
\newtheorem{rem}[thm]{Remark}
\newtheorem*{basic}{Basic Properties}
\newtheorem*{notation}{Notation}
\newtheorem*{cl}{Claim}
\newtheorem{cln}{Claim}

\newcommand{\mC}{{\mathbb C}}

\newcommand{\mN}{{\mathbb N}}

\newcommand{\mQ}{{\mathbb Q}}
\newcommand{\mR}{{\mathbb R}}

\newcommand{\mZ}{{\mathbb Z}}

\newcommand{\Ker}{\mathrm{Ker}\,}

\newcommand{\Spec}{\mathrm{Spec}\,}

\newcommand{\Img}{\mathrm{Im}\,}

\newcommand{\Div}{\mathrm{div}}
\renewcommand{\labelenumi}{(\arabic{enumi})}
\numberwithin{equation}{section}

\begin{document}

\title{An interpretation of multiplier ideals \\ via tight closure}
\author{Shunsuke Takagi}
\address{Graduate School of Mathematical Sciences, University of Tokyo, 3-8-1, Komaba, Meguro, Tokyo 153-8914, Japan}
\email{stakagi@ms.u-tokyo.ac.jp}
\baselineskip 15pt
\footskip = 32pt
\subjclass{Primary 13A35, 14B05; Secondary 14B15, 14E15}
\keywords{Multiplier ideal, Kawamata log terminal, tight closure, strongly F-regular, F-singularities of pairs, modulo $p$ reduction}

\begin{abstract}
Hara \cite{Ha3} and Smith \cite{Sm2} independently proved that in a normal $\mQ$-Gorenstein ring of characteristic $p \gg 0$, the test ideal coincides with the multiplier ideal associated to the trivial divisor. We extend this result for a pair $(R, \Delta)$ of a normal ring $R$ and an effective $\mQ$-Weil divisor $\Delta$ on $\Spec R$. As a corollary, we obtain the equivalence of strongly F-regular pairs and klt pairs.
\end{abstract}

\maketitle
\markboth{SHUNSUKE TAKAGI}{AN INTERPRETATION OF MULTIPLIER IDEALS VIA TIGHT CLOSURE}

\section*{Introduction}
Recently it turned out that there exists a relation between multiplier ideals and tight closure. Precisely speaking, it was proved that some algebraic statements established by multiplier ideals could also be understood via tight closure, for example, Brian\c{c}on-Skoda theorem (see \cite{BS}, \cite{HH1} and \cite{La}), the problem concerning the growth of symbolic powers of ideals in regular local rings (see \cite{ELS} and \cite{HH3}), etc. The purpose of this paper is to give an interpretation of multiplier ideals via tight closure.

The theory of tight closure was introduced by Hochster and Huneke \cite{HH1}, using the Frobenius map in characteristic $p>0$. In this theory, test ideals play a central role.
On the other hand, multiplier ideals, for which we have the strong vanishing theorem, are fundamental tools in birational geometry. Hara \cite{Ha3} and Smith \cite{Sm2} independently proved that in a normal $\mQ$-Gorenstein ring of characteristic $p \gg 0$, the test ideal coincides with the multiplier ideal associated to the trivial divisor. Since the real worth of multiplier ideals is displayed in considering pairs, we attempt to extend this result for pairs. Here, by a \textit{pair}, we mean a pair $(R,\Delta)$ of a normal ring $R$ and an effective $\mQ$-Weil divisor $\Delta$ on $\Spec R$.

Hara generalized the notion of tight closure to that for pairs, which is called $\Delta$-tight closure (cf. \cite[Problem 5.3.2]{HW}). Using the $\Delta$-tight closure operation, we introduce the ideal $\tau(R,\Delta)$ associated to a pair $(R,\Delta)$, which is a generalization of the notion of test ideals. We denote by $\mathcal{J}(Y, \Delta)$ the multiplier ideal of $(Y,\Delta)$. Then, our main result is stated as follows:

\renewcommand{\themainthm}{3.2}
\begin{mainthm}
Let $(R, \mathfrak{m})$ be a normal local ring essentially of finite type over a field of characteristic zero and $\Delta$ an effective $\mQ$-Weil divisor on $Y=\Spec R$ such that $K_Y+ \Delta$ is $\mQ$-Cartier. Then, in characteristic $p \gg 0$,
$$\tau (R, \Delta) = \mathcal{J}(Y, \Delta).$$
\end{mainthm}

When $\Delta=0$, our main theorem coincides with Hara and Smith's result, and they used the strategy which is to reduce to the case where the ring $R$ is quasi-Gorenstein by passing to an index one cover. However we cannot use this strategy, because in general $K_Y$ could be a non-$\mQ$-Cartier divisor, even if $K_Y+\Delta$ is $\mQ$-Cartier. So we give the direct proof which does not use an index one cover.

As a corollary of the main theorem, we get the equivalence of ``F-singularities of pairs'' and ``singularities of pairs.''

The notions of F-regular and F-pure rings, which are closely related to the theory of tight closure, were defined by Hochster and Huneke \cite{HH1} and Hochster and Roberts \cite{HR} respectively. Recently it became clear that F-singularities (F-regular and F-pure rings) correspond to singularities arising in birational geometry (Kawamata log terminal and log canonical singularities). See \cite{Ha2}, \cite{HW}, \cite{MS} and \cite{Sm1}. The notions of Kawamata log terminal (or klt for short) and log canonical (or lc for short) singularities are defined not only for normal rings but also for pairs, and it is these ``singularities of pairs'' that play a very important role in birational geometry.
Therefore Hara and K.-i.~Watanabe \cite{HW} generalized the notions of F-singularities to those for pairs, and they conjectured the equivalence of ``F-singularities of pairs'' and ``singularities of pairs.'' We prove their conjectures.

The ideal $\tau(R,\Delta)$ defines the locus of non-F-regular points of $(R,\Delta)$ in $\Spec R$ under the assumption that $K_R+\Delta$ is a $\mQ$-Cartier divisor. Likewise, the multiplier ideal $\mathcal{J}(\Spec R,\Delta)$ defines the locus of non-klt points of $(\Spec R,\Delta)$. Hence we obtain the following result as a direct consequence of the main result.
\renewcommand{\themaincor}{3.4}
\begin{maincor}[\textup{\cite[Conjecture 5.1.1]{HW}}]
Let $(R, \mathfrak{m})$ be a normal local ring essentially of finite type over a field of characteristic zero and $\Delta$ an effective $\mQ$-Weil divisor on $Y=\Spec R$ such that $K_Y+ \Delta$ is $\mQ$-Cartier. Then, $(Y, \Delta)$ is klt if and only if $(R, \Delta)$ is of strongly F-regular type.
\end{maincor}
\textbf{Acknowledgements.} 
The author wishes to thank Professor Toshiyuki Katsura, his research supervisor, for warm encouragement. He is also grateful to Professor Kei-ichi Watanabe for various comments and many suggestions, Professor Nobuo Hara for helpful advices and valuable information about the subadditivity theorem for the ideal ``$\tau(R,\mathfrak{a})$,'' and Yasunari Nagai for many discussions about birational geometry.

\section{Preliminaries}
\subsection{F-singularities of pairs}
First we briefly review definitions and basic properties on ``F-singularities of pairs,'' which were introduced by Hara and K.-i.~Watanabe. Refer to \cite{HW} for details.

Throughout this paper, all rings are commutative Noetherian integral domains with identity. Let $R$ be an integral domain of characteristic $p > 0$ and $F:R \to R$ the Frobenius map which sends $x$ to $x^p$. 
Since $R$ is reduced, we can identify $F:R \to R$ with the natural inclusion map $R \hookrightarrow R^{1/p}$. $R$ is called \textit{F-finite} if $R \hookrightarrow R^{1/p}$ is a finite map. 
For example, any algebra essentially of finite type over a perfect field is F-finite. We also remark that if $R$ is F-finite, then $R$ is excellent \cite{Ku}.

\begin{notation}
Let $R$ be a normal domain with quotient field $K$.
A $\mQ$-Weil divisor $D$ on $Y=\Spec R$ is a linear combination $D=\sum_{i=1}^{r}a_iD_i$ of irreducible reduced subschemes $D_i \subset Y$ of codimension one with rational coefficients $a_i$. 
The round-up and round-down of $D$ is defined by $\lceil D \rceil=\sum_{i=1}^{r}\lceil a_i \rceil D_i$ and $\lfloor D \rfloor = \sum_{i=1}^{r} \lfloor a_i \rfloor D_i$. We also denote 
$$R(D)=\{x \in K \mid \Div_R(x)+D \geq 0 \}.$$
\end{notation}

\begin{defn}[\textup{\cite[Definition 2.1]{HW}}]
Let $R$ be an F-finite normal domain of characteristic $p>0$ and $\Delta$ an effective $\mQ$-Weil divisor on $\Spec R$.
\begin{enumerate}
\item $(R , \Delta)$ is said to be \textit{F-pure} if the inclusion map $R \hookrightarrow R((q-1)\Delta)^{1/q}$ splits as an $R$-module homomorphism for every $q=p^e$.
\item $(R , \Delta)$ is said to be \textit{strongly F-regular} if for every nonzero element $c \in R$, there exists $q=p^e$ such that $c^{1/q}R \hookrightarrow R((q-1)\Delta)^{1/q} $ splits as an $R$-module homomorphism.
\end{enumerate}
\end{defn}

\begin{rem}
\renewcommand{\labelenumi}{(\roman{enumi})}
\begin{enumerate}
\item $R$ is F-pure (resp. strongly F-regular) if and only if $(R,0)$ is F-pure (resp. strongly F-regular). Refer to \cite{HH1}, \cite{HH2} and \cite{HR} for strongly F-regular and F-pure rings. Here we only note that the following implications hold for rings.

$\text{regular} \Rightarrow \text{strongly F-regular} \Rightarrow \text{normal and Cohen-Macaulay}$

$\hspace{3.5cm} \Downarrow$

$\hspace{2.95cm}\text{F-pure}$

\item We can replace $R((q-1)\Delta)^{1/q}$ by $R(\lceil q\Delta \rceil)^{1/q}$ in the above definition of strong F-regularity.
\end{enumerate}
\end{rem}

\begin{basic}[\textup{\cite[Proposition 2.2]{HW}}] \label{prop:basic(F-sing)}
Let $(R,\Delta)$ be as above.
\renewcommand{\labelenumi}{(\roman{enumi})}
\begin{enumerate}
\item Strong F-regularity implies F-purity.
\item $(R,\Delta)$ is strongly F-regular if and only if for every nonzero element $c \in R$, there exists $q'$ such that $c^{1/q}R \hookrightarrow R(q\Delta)^{1/q}$ splits as an $R$-module homomorphism for all $q= p^e\ge q'$. 
\item If $(R,\Delta)$ is strongly F-regular, then $\lfloor \Delta \rfloor = 0$. 
\item If $(R,\Delta)$ is  F-pure, then $\lceil \Delta \rceil$ is reduced. 
\item If $(R,\Delta)$ is F-pure $($resp.\ strongly F-regular$)$, so is $(R,\Delta')$ for every effective $\mQ$-Weil divisor $\Delta' \leq \Delta$.
\end{enumerate}
\end{basic}

\begin{expl}[\textup{cf. \cite[Theorem 2.5]{Fe} and \cite[Corollary 2.7]{HW}}]
Let $R=$\linebreak[4]
$k[[x_1,\dots,x_n]]$ be an $n$-dimensional complete regular local ring over a field $k$ of characteristic $p>0$ and $\Delta =\Div_R(x_1^{d_1}+\dots+x_n^{d_n})$. Assume that $p$ is sufficiently large and let $t_0=\min \{1,\sum_{i=1}^n\frac{1}{d_i}\}$.
Then, $(R,t\Delta)$ is strongly F-regular if and only if $t<t_0$. If $(R,t\Delta)$ is F-pure, then $t \leq t_0$. When $\sum_{i=1}^n\frac{1}{d_i} > 1$, then $(R,t_0\Delta)$ is always F-pure. On the other hand when $\sum_{i=1}^n\frac{1}{d_i} \leq 1$, then $(R,t_0\Delta)$ is F-pure if $p \equiv 1$ mod $d_i$ for every $i=1,\dots,n$.
\end{expl}

The notions of F-regularity and F-purity are also defined for rings of characteristic zero as follows.
\begin{defn}
Let $R$ be a finitely generated normal domain over a field $k$ of characteristic zero and $\Delta$ an effective $\mQ$-Weil divisor on $\Spec R$. The pair $(R, \Delta)$ is said to be of \textit{F-pure type} (resp. \textit{strongly F-regular type}) if there exist a finitely generated $\mZ$-subalgebra $A$ of $k$, a finitely generated normal $A$-algebra $R_A$ and an effective $\mQ$-Weil divisor $\Delta_A$ on $\Spec R_A$, with a flat structure map $A \to R_A$ such that 
\begin{enumerate}
\item $(A \to R_A) \otimes_A k \cong (k \to R)$ and $\Delta_A \otimes_A k \cong \Delta.$
\item $(R_{\kappa}, \Delta_{\kappa})$ is F-pure (resp. strongly F-regular) for every closed point $s$ in a dense open subset of $\Spec A$, where $\kappa=\kappa(s)$ denotes the residue field of $s \in \Spec A$, $R_{\kappa} = R_A \otimes_A \kappa(s)$ and $\Delta_{\kappa} = \Delta_A \otimes_A \kappa(s)$.
\end{enumerate}
$(R, \Delta)$ is said to be of \textit{dense F-pure type} if in the above condition $(2)$ ``dense open'' is replaced by ``dense.''

\end{defn}
\begin{expl}
Let $R=k[x_1,\dots,x_n]$ be an $n$-dimensional polynomial ring over a field $k$ of characteristic zero and $\Delta =\Div_R(x_1^{d_1}+\dots+x_n^{d_n})$.
Then, $(R,t\Delta)$ is of strongly F-regular type (resp. of dense F-pure type) if and only if $\min \{1,\sum_{i=1}^n\frac{1}{d_i}\} > t$ (resp. $\min \{1,\sum_{i=1}^n\frac{1}{d_i}\} \geq t$).
\end{expl}

\subsection{Birational Geometry}
We recall definitions and fundamental properties of singularities which appear in the Mori theory, and of multiplier ideal sheaves. Refer to \cite{KM} and \cite{La} for details.

Let $Y$ be a normal variety over a field of characteristic zero and $\Delta$ a $\mQ$-Weil divisor on $Y$ such that $K_Y+\Delta$ is $\mQ$-Cartier, that is, $r(K_Y+\Delta)$ is a Cartier divisor for some positive integer $r$, where $K_Y$ is the canonical divisor of $Y$.
Let $f:X \to Y$ be a resolution of singularities such that $ \cup_{i=1}^s E_i+f^{-1}_{*}\Delta$ has simple normal crossing support, where $\mathrm{Exc}(f)=\cup_{i=1}^s E_i$ is the exceptional divisor of $f$ and $f^{-1}_{*}\Delta$ is the strict transform of $\Delta$ in $X$. We denote by $K_X$ the canonical divisor of $X$. Then, for some integers $b_1, \dots ,b_s$,
$$r(K_X+f^{-1}_{*}\Delta) \underset{\text{lin.}}{\sim} f^*(r(K_Y+\Delta)) + \sum_{i=1}^{s}b_iE_i.$$
Hence we have
$$K_X+f^{-1}_{*}\Delta \underset{\text{$\mQ$-lin.}}{\sim} f^*(K_Y+\Delta) + \sum_{i=1}^{s}a_iE_i,$$
where $a_i = b_i/r$ ($i=1,\dots, s$).
\begin{defn}
Under the same notation as above:
\begin{enumerate}
\item We say that the pair $(Y,\Delta)$ is \textit{Kawamata log terminal} (or \textit{klt} for short) if $a_i > -1$ for every $i=1,\dots, s$ and $\lfloor \Delta \rfloor \leq 0$.
\item We say that the pair $(Y,\Delta)$ is \textit{log canonical} (or \textit{lc} for short) if $a_i \geq -1$ for every $i=1,\dots, s$ and the coefficient of $\Delta$ in each irreducible component is less than or equal to one.
\item The \textit{multiplier ideal sheaf} $\mathcal{J}(Y, \Delta)$ associated to $\Delta$ is defined to be
$$\mathcal{J}(Y, \Delta)=f_*\mathcal{O}_X(\lceil K_X -f^*(K_Y+\Delta) \rceil).$$
\end{enumerate}
\end{defn}
\begin{rem}
\renewcommand{\labelenumi}{(\roman{enumi})}
\begin{enumerate}
\item The above definitions do not depend on the choice of a desingularization $f:X \to Y$.
\item When $\Delta$ is effective, $\mathcal{J}(Y, \Delta)$ is indeed an ideal sheaf. However in case $\Delta$ is not effective, it is generally not a submodule of $\mathcal{O}_Y$ but a fractional ideal sheaf.
\end{enumerate}
\end{rem}
\begin{basic}In the situation of the above definition:
\renewcommand{\labelenumi}{(\roman{enumi})}
\begin{enumerate}
\item For every $\mQ$-Weil divisor $\Delta' \leq \Delta$ on $Y$, $\mathcal{J}(Y,\Delta') \supseteq \mathcal{J}(Y,\Delta)$.
\item For every Cartier divisor $\Delta'$ on $Y$, $\mathcal{J}(Y,\Delta + \Delta') = \mathcal{J}(Y,\Delta) \otimes_Y \mathcal{O}_Y (- \Delta')$.
\item The pair $(Y, \Delta)$ is klt if and only if $\mathcal{J}(Y,\Delta) \supseteq \mathcal{O}_Y$. In particular when $\Delta$ is effective, $(Y, \Delta)$ is klt if and only if $\mathcal{J}(Y,\Delta)=\mathcal{O}_Y$.
\end{enumerate}
\end{basic}

\begin{prop}[\textup{cf. \cite[Proposition 5.20]{KM}, \cite[Example 5.26]{La}}]\label{etale in zero}
Let $f:X \to Y$ be a finite covering of normal varieties which is \'{e}tale in codimension one and $\Delta$ an effective $\mQ$-Weil divisor on $Y$. Then
$$\mathcal{J}(Y,\Delta)=\mathcal{J}(X,f^{*}\Delta) \cap \mathcal{O}_Y.$$
\end{prop}

\begin{prop}[\textup{\cite{DEL}}]\label{DEL}
\begin{enumerate}
\item \textup{(Restriction Theorem)} Let $Y$ be a normal variety, $\Delta$ an effective divisor on Y such that $K_Y+\Delta$ is $\mQ$-Cartier, and $H$ a reduced Cartier divisor which is not in the support of $\Delta$. Assume that $H$ is a normal variety. Then $$\mathcal{J}(H,\Delta|_H) \subseteq \mathcal{J}(Y,\Delta) \cdot \mathcal{O}_H.$$
\item \textup{(Subadditivity Theorem)} Let $Y$ be a non-singular quasi-projective variety, and $\Delta_1$ and $\Delta_2$ be any two effective $\mQ$-divisors on $Y$. Then
$$\mathcal{J}(Y,\Delta_1+\Delta_2) \subseteq \mathcal{J}(Y,\Delta_1) \cdot \mathcal{J}(Y,\Delta_2).$$
\end{enumerate}
\end{prop}

\begin{expl}[\textup{\cite{La}}]
\begin{enumerate}
\item When $Y$ is a non-singular variety and $\Delta$ is a $\mQ$-Weil divisor on $Y$ with simple normal crossing support, then $\mathcal{J}(Y,\Delta) = \mathcal{O}_Y(-\lfloor \Delta \rfloor)$.
\item Let $Y={\mC}^n$ with coordinates $x_1,\dots,x_n$ and $\Delta = \Div_Y(x_1^{d_1}+\dots+x_n^{d_n})$. Then, $(Y,t\Delta)$ is klt (resp. lc) if and only if $\min \{1,\sum_{i=1}^n\frac{1}{d_i}\} >t$ (resp. $\min \{1,\sum_{i=1}^n\frac{1}{d_i}\} \geq t$).
\item Let $Y={\mC}^d$ with coordinates $x_1,\dots,x_d$ and $\Delta=\Div_Y(x_1^{d+1}+\dots+x_d^{d+1})$. Then, $\mathcal{J}(Y, \frac{d}{d+1}\Delta) = (x_1,\dots,x_d)$.
\end{enumerate}
\end{expl}

\section{$\Delta$-tight closure}
In this section, we introduce the notion of $\Delta$-tight closure which is suggested in \cite[Problem 5.3.2]{HW} and see that the $\Delta$-tight closure operation satisfies properties similar to those of the ``usual'' tight closure operation (see \cite{Hu} for the theory of ``usual'' tight closure). 
Moreover, using the $\Delta$-tight closure operation, we define the ideal $\tau(R,\Delta)$, which is a generalization of the ``usual'' test ideal.

\begin{notation}
Let $R$ be a normal domain of characteristic $p>0$ and $\Delta$ an effective $\mQ$-Weil divisor on $\Spec R$. 
\begin{itemize}
\item We always use the letter $q$ (resp. $q'$, $q''$) for a power $p^e$ (resp. $p^{e'}$, $p^{e''}$) of $p$. 
\item For any ideal $I$ in $R$, we denote by $I^{[q]}$ the ideal of $R$ generated by the $q$ th powers of elements of $I$.
\item For any divisorial ideal $J$ of R (i.e., $J=R(D)$ for some unique integral Weil divisor $D$), we denote by $J^{(m)}$ the reflexive hull of $J^m$. If $J=R(D)=R(\lfloor D \rfloor)$, then $J^{(m)}=R(m\lfloor D \rfloor)$. 
\item The notation ${}^eR((q-1)\Delta)$ denotes $R((q-1)\Delta)$ itself, but viewed as an $R$-module via the $e$-times Frobenius map $F^e:R \to R((q-1)\Delta)$.\item When $(R,\mathfrak{m})$ is local, we denote by $E_R$ the injective hull of the residue field $R/\mathfrak{m}$. 
\end{itemize}
\end{notation}

\begin{defn}[\textup{cf. \cite[Problem 5.3.2]{HW}}]
Let $N \subseteq M$ be modules over an F-finite normal domain $R$ of characteristic $p>0$ and $\Delta$ an effective $\mQ$-Weil divisor on $\Spec R$.
We denote by $F^e:M=M \otimes _R R \to M \otimes _R {}^eR((q-1)\Delta)$ the $e$-times Frobenius map induced on $M$ which sends $z \in M$ to $z^q:=F^e(z)=z \otimes 1 \in  M \otimes _R {}^eR((q-1)\Delta)$. 
Set $N^{[q]\Delta}_M := \Img (F^e(N) \to F^e(M))$. 
Then the \textit{$\Delta$-tight closure} $N_M^{*\Delta} \subseteq M$ of $N$ in $M$ is defined as follows: $z \in N_M^{*\Delta}$ if and only if there exists a nonzero element $c \in R$ such that $cz^q:= z \otimes c \in N^{[q]\Delta}_M$ for all $q=p^e \gg 0$.
$$\xymatrix{
N=N \otimes _R R \ar@{^{(}->}[r] \ar[d]^{F^e} & M=M \otimes _R R \ar[d]^{F^e}\\
N \otimes _R {}^eR((q-1)\Delta) \ar[r] & M \otimes _R {}^eR((q-1)\Delta)
}$$
Moreover, the \textit{finitistic $\Delta$-tight closure} $N_M^{*\Delta fg} \subseteq M$ of $N$ in $M$ is defined to be $N_M^{*\Delta fg} := \underset{M'}{\cup} N^{*\Delta}_{M'}$, where $M'$ runs through all finitely generated $R$-submodules of $M$ which contain $N$.
\end{defn}

\begin{rem}
\renewcommand{\labelenumi}{(\roman{enumi})}
\begin{enumerate}
\item When $\Delta=0$, $\Delta$-tight closure coincides with ``usual'' tight closure. 
\item In general, $N^{*\Delta}_M \subsetneq (N^{*\Delta}_M)^{*\Delta}_M$. In this sense, the $\Delta$-tight closure operation is not a ``closure operation.''
\item Let $I$ be an ideal in $R$. Then, $I^{[q]\Delta}_R=I^{[q]}R((q-1)\Delta)$.
\item We can replace ${}^e R((q-1)\Delta)$ by ${}^e R(\lceil q\Delta \rceil)$ in the above definition.
\item In general, $N_M^{*\Delta fg} \subseteq N_M^{*\Delta}$. If $M$ itself is finitely generated, then $N_M^{*\Delta fg}=N_M^{*\Delta}$.
\end{enumerate}
\end{rem}

\begin{basic}
In the situation of the above definition,
\renewcommand{\labelenumi}{(\roman{enumi})}
\begin{enumerate}
\item $N \subseteq N^{*\Delta}_M$.
\item $N^{*\Delta}_M/N = 0^{*\Delta}_{M/N}$.
\item For any effective $\mQ$-Weil divisor $\Delta' \leq \Delta$ on $\Spec R$, $N^{*\Delta'}_M \subseteq N^{*\Delta}_M$.
\item For any effective Cartier divisor $\Delta'$ on $\Spec R$, $N^{*(\Delta+\Delta')}_M = N^{*\Delta}_M : R( -\Delta' )$.
\item If $(R,\Delta)$ is strongly F-regular, then $I^{*\Delta}=I$ for every ideal $I$ in $R$.
\end{enumerate}
\end{basic}

Strong F-regularity can be characterized via $\Delta$-tight closure. 
\begin{lem}\label{lem:0_E^{Delta'*}=0}
Let $(R, \mathfrak{m})$ be an F-finite normal local ring of characteristic $p>0$ and $\Delta$ an effective $\mQ$-Weil divisor on $\Spec R$. Then $(R,\Delta)$ is strongly F-regular if and only if $ 0_E^{*\Delta}=0$, where $E=E(R/\mathfrak{m})$ is the injective hull of the residue field $R/\mathfrak{m}$.
\end{lem}
\begin{proof}
The proof is essentially the same as that for the no boundary case \cite[Proposition 2.1]{Ha1}.

Assume that $(R,\Delta)$ is strongly F-regular and fix any element $z \in 0_E^{*\Delta}$.
Then there exists a nonzero element $c \in R$ such that $cz^q:=cF^e(z)=0$ for all $q=p^e \gg 0$.
Let
$$ \phi_c^{(e)}:\mathrm{Hom}_R(R((q-1)\Delta)^{1/q}, R) \to \mathrm{Hom}_R(R,R)=R $$
be an $R$-module homomorphism induced by the $R$-linear map $R \xrightarrow{c^{1/q}} R((q-1)\Delta)^{1/q}$ for each $q=p^e$. 
Since $(R,\Delta)$ is strongly F-regular, $\phi_c^{(e)}$ is surjective for all $q = p^e \gg 0$. 
Since the $R$-module homomorphism $cF^e:E \to E \otimes_R {}^eR((q-1)\Delta)$ which sends $z$ to $cz^q$ is the Matlis dual of $\phi_c^{(e)}$, $cF^e$ is injective for every $q = p^e \gg 0$. Hence we have $z=0$. 

Conversely, suppose that $ 0_E^{*\Delta}=0$, and fix any nonzero element $c \in R$. If $z$ is a nonzero element of the socle $(0:\mathfrak{m})_E$ of $E$, then there exists $q=p^e$ such that $cF^e(z) \ne 0$. 
Since $(0:\mathfrak{m})_E$ is an one-dimensional $R/\mathfrak{m}$-vector space, we can take $q$ which works for every $z \in (0:\mathfrak{m})_E$. Then $cF^e$ is injective on $(0:\mathfrak{m})_E$. 
Since $E$ is an essential extension of $(0:\mathfrak{m})_E$, $cF^e$ itself is injective. 
Taking the Matlis dual of $cF^e$, we know that $\phi_c^{(e)}$ is surjective, namely $(R,\Delta)$ is strongly F-regular.
\end{proof}

We introduce $\Delta$-test elements which are very useful to show the propositions about $\Delta$-tight closure. When $\Delta$ is zero, a $\Delta$-test element is nothing but a test element in the theory of ``usual'' tight closure. 
\begin{defn}
 Let $R$ be an F-finite normal domain of characteristic $p>0$ and $\Delta$ an effective $\mQ$-Weil divisor on $\Spec R$.
A nonzero element $c \in R$ is called a \textit{$\Delta$-test element} if for each ideal $I$ in $R$, $x \in I^{*\Delta}$ if and only if $cx^q \in I^{[q]}R(\lceil q\Delta \rceil)$ for all $q=p^e$.
\end{defn}

By the following lemma, which is a generalization of \cite[Theorem 3.3]{HH2}, $\Delta$-test elements always exist.
\begin{lem}\label{lem:test-element}
Let $R$ be an F-finite normal domain of characteristic $p>0$.
Let $c \in R$ be any nonzero element such that the localization $R_c$ with respect to $c$ is strongly F-regular.
\begin{enumerate}
\renewcommand{\labelenumi}{(\arabic{enumi})}
\item For every effective divisor $\Delta$ on $R$, $(R , \Delta)$ is strongly F-regular if and only if there exists $q=p^e$ such that $c^{1/q}R \hookrightarrow R(\lceil q\Delta \rceil)^{1/q}$ splits as an $R$-module homomorphism.
\item For every effective divisor $\Delta$ on $R$, $c^n$ is a $\Delta$-test element for some positive integer $n$.
\end{enumerate}
\end{lem}
\begin{proof}
Let $d \in R$ be any nonzero element.
\setcounter{cln}{0}
\begin{cln}
For some $q'$ and sufficiently large $q''$, there exists an $R^{1/qq''}$-module homomorphism
$$ R(\lceil qq'q''\Delta \rceil)^{1/qq'q''} \to R(\lceil qq''\Delta \rceil)^{1/qq''}, \quad d^{1/qq'q''} \mapsto c^{1/q}.$$
\end{cln}
\begin{proof}[Proof of Claim $1$]
Since $R_c$ is strongly F-regular, by the proof of \cite[Theorem 3.3]{HH2}, for some $q'$ and sufficiently large $q''$, there exists an $R$-module homomorphism
$$R^{1/q'} \to R, \quad d^{1/q'} \mapsto c^{q''}.$$
Tensoring with $R(\lceil qq''\Delta \rceil)$ and taking $qq''$-th roots, for each $q=p^e$, we get an $R^{1/qq''}$-linear map
$$ R(\lceil qq'q''\Delta \rceil)^{1/qq'q''} \to R(\lceil qq''\Delta \rceil)^{1/qq''}, \quad d^{1/qq'q''} \mapsto c^{q''/qq''}=c^{1/q}.$$
\end{proof}

First we will show $(1)$. Suppose that there exists $q=p^e$ such that $c^{1/q}R \hookrightarrow R(\lceil q\Delta \rceil)^{1/q}$ splits as an $R$-module homomorphism, that is, there exists an $R$-linear map $R(\lceil q\Delta \rceil)^{1/q} \to R$ sending $c^{1/q}$ to $1$. 
We may replace $q''$ in Claim $1$ suitably (for example, $q''=q^n$ for some positive integer $n$), and then there exists an $R$-module homomorphism $R(\lceil q''\Delta \rceil)^{1/{q''}} \to R$ which sends $1$ to $1$. 
Tensoring this with $R(\lceil q\Delta \rceil)$ and taking $q$-th roots, we obtain an $R^{1/q}$-linear map $R(\lceil qq''\Delta \rceil)^{1/qq''} \to R(\lceil q\Delta \rceil)^{1/q}$ which sends $1$ to $1$. 
By composing these maps with the $R^{1/qq''}$-linear map in Claim $1$, we get the following $R$-module homomorphism
$$R(\lceil qq'q''\Delta \rceil)^{1/qq'q''} \to R(\lceil qq''\Delta \rceil)^{1/qq''}  \to R(\lceil q \Delta \rceil)^{1/q}  \to R, $$
$$d^{1/qq'q''}  \longmapsto  c^{1/q}  \longmapsto c^{1/q}  \longmapsto 1.$$
This establishes $(1)$.

Next we will prove $(2)$. 
When $d=1$, then we can take $p$ as $q'$ in the proof of Claim $1$, namely, there exists an $R$-module homomorphism
$$h:R^{1/p} \to R, \quad 1 \mapsto c^{q''}.$$
We replace $c^{q''}$ by $c$, and then it is enough to show that $c^3$ is a $\Delta$-test element.
\begin{cln}
For every $q=p^e$, there exists an $R$-module homomorphism
$$g_e:R(q\lceil \Delta \rceil)^{1/q} \to R(\lceil \Delta \rceil), \quad 1 \mapsto c^2.$$
\end{cln}
\begin{proof}[Proof of Claim $2$]
When $q=p$, then tensoring $h$ with $R(\lceil \Delta \rceil)$, we obtain an $R$-linear map
$$h_1:R(p\lceil \Delta \rceil)^{1/p} \to R(\lceil \Delta \rceil) \quad 1 \mapsto c.$$
Hence we may set $g_1=c \cdot h_1$. 
Suppose that the assertion holds for $q=p^e$. 
Then, by tensoring with $R((p-1)\lceil \Delta \rceil)$ and taking $p$-th roots, we obtain an $R^{1/p}$-linear map $R(pq\lceil \Delta \rceil)^{1/pq} \to R(p\lceil \Delta \rceil)^{1/p}$ which sends $1$ to $c^{2/p}$. 
We may compose this with an $R^{1/p}$-module homomorphism $R(p\lceil \Delta \rceil)^{1/p} \to R(p\lceil \Delta \rceil)^{1/p}$ which sends $1$ to $c^{(p-2)/p}$, and then with $h_1$.
$$R(pq\lceil \Delta \rceil)^{1/pq} \to R(p\lceil \Delta \rceil)^{1/p} \to R(p \lceil \Delta \rceil)^{1/p} \to R,$$
$$1 \longmapsto c^{2/p} \longmapsto c \longmapsto c^2.$$
This is a required homomorphism for $pq=p^{e+1}$.
\end{proof}

Thanks to Claim $2$, tensoring $g_{e''}$ with $R(\lceil q\Delta \rceil-\lceil \Delta \rceil)$ and taking $q$-th roots, we obtain an $R^{1/q}$-module homomorphism
$$R(\lceil qq''\Delta \rceil)^{1/qq''} \to R(\lceil q\Delta \rceil)^{1/q}, \quad 1 \mapsto c^{2/q}.$$
Therefore, by Claim $1$, there exists an $R^{1/q}$-linear map
$$R(\lceil qq'q''\Delta \rceil)^{1/qq'q''} \to R(\lceil qq''\Delta \rceil)^{1/qq''} \to R(\lceil q\Delta \rceil)^{1/q},$$
$$d^{1/qq'q''} \longmapsto c^{1/q} \longmapsto c^{3/q}.$$
Thus, $dx^{qq'q''} \in$ $I^{[qq'q'']}R(\lceil qq'q''\Delta \rceil)$ implies $c^3x^q \in I^{[q]}R(\lceil q \Delta \rceil)$. It follows that $c^3$ is a $\Delta$-test element.
\end{proof}

Now we define the ideal $\tau(R,\Delta)$.
\begin{defn}
Let $R$ be an F-finite normal domain of characteristic $p>0$ and $\Delta$ an effective $\mQ$-Weil divisor on $Y=\Spec R$. 
Then the ideal $\tau(R, \Delta)$ is defined to be $\tau(R, \Delta):=\bigcap_M \mathrm{Ann}_{R}(0_{M}^{*\Delta}) \subseteq R$, where $M$ runs through all finitely generated $R$-modules.
\end{defn}

\begin{rem}
When $\Delta=0$, $\tau(R,\Delta)$ coincides with the ``usual'' test ideal which is generated by ($\Delta$-)test elements. 
However even if $K_Y+\Delta$ is $\mQ$-Cartier, $\tau(R,\Delta)$ may not be generated by $\Delta$-test elements. Therefore We do not call $\tau(R,\Delta)$ the $\Delta$-test ideal.
\end{rem}

\begin{basic}
In the situation of the above definition:
\renewcommand{\labelenumi}{(\roman{enumi})}
\begin{enumerate}
\item For any effective $\mQ$-Weil divisor $\Delta' \leq \Delta$ on $\Spec R$, we have $\tau(R,\Delta') \supseteq \tau(R,\Delta)$.
\item For any effective Cartier divisor $\Delta'$ on $\Spec R$, we have $\tau(R,\Delta + \Delta') = \tau(R,\Delta) \otimes_R R(-\Delta')$.
\item If $(R,\Delta)$ is strongly F-regular, then $\tau(R,\Delta)=R$.
\end{enumerate}
\end{basic}

\begin{thm}\label{lem:test-ideal}
Let $R$ be an F-finite normal domain of characteristic $p>0$ and $\Delta$ an effective $\mQ$-Weil divisor on $Y=\Spec R$. Let $E=\bigoplus_{\mathfrak{m}}E(R/\mathfrak{m})$, where $\mathfrak{m}$ runs through all maximal ideals of $R$.
\begin{enumerate}
\item $$\tau(R, \Delta) = \bigcap_I (I:I^{*\Delta}) =\mathrm{Ann}_R(0_E^{*\Delta fg}),$$
where the intersection in the middle term is taken over all ideals $I$ in $R$.
\item  If $K_Y + \Delta$ is a $\mQ$-Cartier divisor, then
$$\tau(R, \Delta) =\mathrm{Ann}_R(0_E^{*\Delta}).$$
\end{enumerate}
\end{thm}
\begin{proof}
By the same argument as in the proof of \cite[Proposition 8.23]{HH1}, we see that $\bigcap_M \mathrm{Ann}_R(0_M^{*\Delta}) =\bigcap_I (I:I^{*\Delta})=\mathrm{Ann}_R (0_E^{*\Delta fg})$. So we will prove $(2)$. We may assume that $(R,\mathfrak{m})$ is local. Since $0_E^{*\Delta} \supseteq 0_E^{*\Delta fg}$, it is enough to show that $\mathrm{Ann}_R(0_E^{*\Delta}) \supseteq \bigcap_I (I:I^{*\Delta})$ under the assumption that $K_Y + \Delta$ is $\mQ$-Cartier.

We use the same strategy as that of \cite{Ha3}. For a sequence of elements $\mathbf{x}=x_1,\dots,x_d$ of $R$ and a positive integer $t$, we write $\mathbf{x}^t=x_1^t,\dots,x_d^t$. For an $R$-module $M$, we denote
$$\mathcal{K}(\mathbf{x},t,M):=\Ker \left(\frac{M}{(\mathbf{x})M} \xrightarrow{(x_1 \cdots x_d)^{t-1}} \frac{M}{(\mathbf{x}^t)M} \right).$$
We also denote
$$\mathcal{K}(\mathbf{x},\infty,M):= \bigcup_{t \in \mN} \mathcal{K}(\mathbf{x},t,M).$$

\begin{cl}
Let $(R,\mathfrak{m})$ be a $d$-dimensional normal local ring of characteristic $p > 0$ and $J \subseteq R$ a divisorial ideal. Let $\mathbf{x}=x_1, \dots, x_d$ be a system of parameters of $R$. Suppose that there exist a nonzero element $c \in R$ and integer $t_0 \geq 2$ such that
$$c \mathcal{K}(\mathbf{x}^{qs},\infty, J^{[q]}R((q-1)\Delta)) \subseteq \mathcal{K}(\mathbf{x}^{qs},t_0, J^{[q]}R((q-1)\Delta))$$
for all $s \geq 1$ and $q=p^e \gg 0$.
Then $\mathrm{Ann}_R(0_{H^d_{\mathfrak{m}}(J)}^{*\Delta}) \supseteq \bigcap_I (I:I^{*\Delta})$.
\end{cl}
\begin{proof}[Proof of Claim]
Let $\xi \in 0_{H^d_{\mathfrak{m}}(J)}^{*\Delta}$, i.e., there exists a nonzero element $d \in R$ such that $d\xi^q=0 \in H^d_{\mathfrak{m}}(J) \otimes_R {}^eR((q-1)\Delta)$ for every $q=p^e \gg 0$. Note that $H^d_{\mathfrak{m}}(J) \cong \varinjlim J/(\mathbf{x}^{t})J$, where the direct limit map $J/(\mathbf{x}^{t})J \to J/(\mathbf{x}^{t+1})J$ is the multiplication by $x_1x_2 \cdots x_d$. Therefore $\xi$ is represented by $z$ mod $(\mathbf{x}^{s})J \in J/(\mathbf{x}^{s})J$ for some $z \in J$ and $s \ge 1$. We may assume that $s=1$.

Now the natural map 
$$\frac{J}{(\mathbf{x}^t)J} \otimes_R {}^eR((q-1)\Delta) \to \frac{J^{[q]}R((q-1)\Delta)}{(\mathbf{x}^{qt})J^{[q]}R((q-1)\Delta)}$$ 
induces a map
$$H^d_{\mathfrak{m}}(J) \otimes_R {}^eR((q-1)\Delta) \to \varinjlim \frac{J^{[q]}R((q-1)\Delta)}{(\mathbf{x}^{qt})J^{[q]}R((q-1)\Delta)}$$
which sends $\xi^q$ to the class of $z^q$ mod $(\mathbf{x}^{q})J^{[q]}R((q-1)\Delta)$. 
Thus, for every $q=p^e \gg0$, 
$$\textup{ class of }dz^q\textup{ mod }(\mathbf{x}^{q})J^{[q]}R((q-1)\Delta)=0 \in \varinjlim \frac{J^{[q]}R((q-1)\Delta)}{(\mathbf{x}^{qt})J^{[q]}R((q-1)\Delta)},$$
whence $dz^q \in \mathcal{K}(\mathbf{x}^q,\infty,J^{[q]}R((q-1)\Delta))$. By our assumption, we have $cdz^q \in \mathcal{K}(\mathbf{x}^{q},t_0, J^{[q]}R((q-1)\Delta))$, namely, 
$$cdz^q(x_1 \cdots x_d)^{q(t_0-1)} \in (\mathbf{x}^{qt_0})J^{[q]}R((q-1)\Delta)$$
for every $q=p^e \gg 0$. It implies that $z(x_1 \cdots x_d)^{(t_0-1)} \in ((\mathbf{x}^{t_0})J)^{*\Delta}$. Fix any element $a \in \bigcap_I (I:I^{*\Delta})$. Then $az(x_1 \cdots x_d)^{(t_0-1)} \in (\mathbf{x}^{t_0})J$. Thus it follows that
$$a\xi=\textup{class of }az\textup{ mod }(\mathbf{x})J =0 \in \varinjlim J/(\mathbf{x}^{t})J.$$
Hence the element $a$ is contained in $\mathrm{Ann}_R((0_{H^d_{\mathfrak{m}}(J)}^{*\Delta})$.
\end{proof}

Let $J=R(D) \subset R$ be a divisorial ideal which is isomorphic to the canonical module $\omega_R$ as an $R$-module. 
Now we will prove that there exist a nonzero element $d \in R$ and a system of parameters $\mathbf{x}=x_1,\dots,x_d$ of $R$ such that
$$d \mathcal{K}(\mathbf{x}^{qs},\infty, J^{[q]}R((q-1)\Delta)) \subseteq \mathcal{K}(\mathbf{x}^{qs},2, J^{[q]}R((q-1)\Delta))$$
for all $s \geq 1$ and $q=p^e \gg 0$.
Since $r(D+\Delta)$ is a Cartier divisor for some positive integer $r$, let $R(r(D+\Delta))=yR$. Fix any $a \in R(-r\Delta)$, and let $x_1:=ay \in J$. Then, by \cite[Lemma 4.3]{Wi}, there exist an element $x_2 \in R$ which is not in any minimal prime divisor of $x_1$ and $c \in J$ such that $x_2^n J^{(n)} \subseteq c^nR$ for all $n > 0$. The sequence $x_1$, $x_2$ can be extended to a system of parameters $\mathbf{x}=x_1, \dots, x_d$ for $R$.
Now given any power $q=p^e$, write $q-1=kr+i$ for integers $k$ and $i$ with $0 \leq i \leq r-1$. Then we have
\begin{align*} 
c^r x_2^q R(krD+(q-1)\Delta) & \subseteq c^r x_2^{kr} (J^{(kr)} \otimes_R R((q-1)\Delta))^{**} \\ 
& \subseteq c^{kr+r}R((q-1)\Delta)\\
& \subseteq J^{[q]}R((q-1)\Delta),
\end{align*}
where $(\,)^{**}$ denotes the reflexive hull.
Since $x_1^q \in J^{[q]}$, this implies that
$$c^{r}(x_1 \cdots x_{i-1}x_{i+1} \cdots x_d)^{qs}R(krD+(q-1)\Delta) \subseteq J^{[q]}R((q-1)\Delta)$$
for every $s \geq 1$ and $i=1,\dots,n$.
Therefore, letting $c_1=c^{r}$, we have
\begin{align*}
c_1 \cdot \Img \bigg( H^{d-1}\Bigl(\mathbf{x}^{qst} \,;\, \frac{R(krD+(q-1)\Delta)}{J^{[q]}R((q-1)\Delta)}\Bigr) \\
&\hspace*{-1in}\to H^{d-1}\Bigl(\mathbf{x}^{qst+qs}\,;\, \frac{R(krD+(q-1)\Delta)}{J^{[q]}R((q-1)\Delta)}\Bigr) \bigg)=0
\end{align*}
for all $s,t \geq 1$.
On the other hand, let $c' \in R$ be a test element (By Lemma \ref{lem:test-element}, such $c'$ always exists). If $z \in R(i\Delta)$ is an element such that $z$ mod $(\mathbf{x}^{qs})R(i\Delta) \in \mathcal{K}(\mathbf{x}^{qs},\infty, R(i\Delta))$, then 
$$az \in (x_1^{qst},\dots,x_n^{qst}):(x_1 \cdots x_n)^{qs(t-1)}$$ 
for some $t \geq 1$, so $az \in (x_1^{qs},\dots,x_n^{qs})^{*}$ by colon-capturing \cite{HH1}. Hence, letting $c_2= a \cdot c'$, we have
$$c_2\mathcal{K}(\mathbf{x}^{qs},\infty, R(krD+(q-1)\Delta)) \cong c_2\mathcal{K}(\mathbf{x}^{qs},\infty, R(i\Delta))=0.$$
Thus, applying \cite[Lemma A.3]{Ha3} to the exact sequence
$$0 \to J^{[q]}R((q-1)\Delta) \to R(krD+(q-1)\Delta) \to \frac{R(krD+(q-1)\Delta)}{J^{[q]}R((q-1)\Delta)} \to 0,$$
we see that 
$$c_1c_2\mathcal{K}(\mathbf{x}^{qs},\infty, J^{[q]}R((q-1)\Delta)) \subseteq \mathcal{K}(\mathbf{x}^{qs},2, J^{[q]}R((q-1)\Delta))$$
for all $s,t \geq 1$. Thanks to the claim, we obtain the assertion.
\end{proof}

\begin{rem}
Theorem \ref{lem:test-ideal} implies that under the condition that $R$ is complete and $K_Y+\Delta$ is $\mQ$-Cartier, the finitistic $\Delta$-tight closure of zero coincides with the $\Delta$-tight closure of zero in $E$. In case $\Delta=0$, this holds without the assumption that $R$ is complete (See \cite{AM}, \cite{Mc} and \cite{Sm2}).
Hence we believe that this coincidence of the finitistic $\Delta$-tight closure and the $\Delta$-tight closure is true, even if $R$ is not necessarily complete.
\end{rem}

By Lemma \ref{lem:0_E^{Delta'*}=0}, if $K_R+\Delta$ is $\mQ$-Cartier, then $\tau(R,\Delta)$ defines the locus of non-F-regular points of $(R,\Delta)$ in $\Spec R$.
\begin{cor}
Let $(R, \mathfrak{m})$ be an F-finite normal local ring of characteristic $p>0$ and $\Delta$ an effective $\mQ$-Weil divisor on $\Spec R$ such that $K_Y + \Delta$ is $\mQ$-Cartier. Then $(R,\Delta)$ is strongly F-regular if and only if $\tau(R,\Delta)=R$, which is equivalent to the condition that $I^{*\Delta}=I$ for every ideal $I$ in $R$.
\end{cor}

From now on, we treat only the case where $K_R+\Delta$ is $\mQ$-Cartier, whence $\tau(R,\Delta)=\mathrm{Ann}_R(0^{*\Delta}_E)$.
The following proposition corresponds to Proposition \ref{etale in zero}.
\begin{prop}\label{etale}
Let $(R,\mathfrak{m}) \hookrightarrow (S,\mathfrak{n})$ be a finite local homomorphism of F-finite normal local rings of characteristic $p >0$ which is \'{e}tale in codimension one. Let $\Delta_R$ be an effective $\mQ$-Weil divisor on $\Spec R$ such that $K_R+\Delta_R$ is $\mQ$-Cartier, and $\Delta_S$ be the pullback of $\Delta_R$ by the induced morphism $\pi:\Spec S \to \Spec R$. Assume that $\mathrm{deg}$ $\pi$ is not divisible by $p$. Then
$$\tau(R,\Delta_R) = \tau(S,\Delta_S) \cap R.$$
\end{prop}
\begin{proof}
Note that, by \cite[Proposition 5.7]{KM}, $R$ is a direct summand of $S$ as an $R$-module. Therefore we consider $E_R$ as a direct summand of $E_S=S \otimes_R E_R$. 
If $\xi \in 0^{*\Delta_R}_{E_R}$, then there exists a nonzero element $c \in R$ such that $c\xi^q =0$ in $E_R \otimes_R {}^eR((q-1)\Delta_R)$ for all $q=p^e \gg 0$, so $c\xi^q =0$ in $E_S \otimes_R {}^eR((q-1)\Delta_R) = E_S \otimes_S {}^eS((q-1)\Delta_S)$ (this equality follows from the assumption that $R \hookrightarrow S$ is \'{e}tale in codimension one). This implies that $\tau(S,\Delta_S) \cap R \subseteq \tau(R,\Delta_R)$.

Conversely, let $c$ be a nonzero element of $\tau(R,\Delta_R)$, and fix any nonzero element $d \in R$. Let $F^e_R:E_R \to E_R \otimes_R {}^eR((q-1)\Delta_R)$ (resp. $F^e_S:E_S \to E_S \otimes_S {}^eS((q-1)\Delta_S)$) be the $e$-times Frobenius map induced on $E_R$ (resp. $E_S$). 
Then $c \cdot \underset{e \geq e'}{\cap} \Ker dF^e_R =0$ for every $q'=p^{e'}$. Since $E_R$ is Artinian, there exists $q''=p^{e''}$ such that $\underset{e \geq e'}{\cap} \Ker dF^e_R = \underset {e'' \geq e \geq e'}{\cap} \Ker dF^e_R$. 
As in (Step $1$) of the proof of Theorem \ref{thm:test subset multi}, for every ${e''} \geq e \geq {e'}$, there exist $c_e \in R$ and an $R$-module homomorphism $ {\phi_e}':R((q-1)\Delta)^{1/q} \to R$ sending $d^{1/q}$ to $c_e$ such that $\underset{{e''} \geq e \geq {e'}}{\sum}c_e=c$.
Since $R \hookrightarrow S$ is \'{e}tale in codimension one, by tensoring this homomorphism with $E_S$ over $R$, we get the following commutative diagram (cf. \cite[Theorem 4.8]{HW})
$$\xymatrix 
{E_S \ar[r]^{\hspace{-1.5cm} dF^e_S} \ar[dr]^{c_e} & E_S \otimes_S {}^e S((q-1)\Delta_S) \ar[d]^{\phi_e} \\
& E_S, \\
}$$
where the $S$-linear map $\phi_e:E_S \otimes_S {}^e S((q-1)\Delta_S) \to E_S$ sends $z \otimes d$ to $c_e z$.
Therefore $\Ker dF^e_S \subseteq (0:c_e)_{E_S}$ for every ${e''} \geq e \geq {e'}$.
Since $\underset{{e''} \geq e \geq {e'}}{\sum}c_e=c$, we have $\underset {e'' \geq e \geq e'}{\cap} \Ker dF^e_S \subseteq (0:c)_{E_S}$, and hence $c \cdot \underset{e \geq e'}{\cap} \Ker dF^e_S =0$ for every $d \in R$ and $q'=p^{e'}$. 
Since $R \hookrightarrow S$ is a finite extension of normal domains, so $dS \cap R \neq 0$ for all nonzero elements $d \in S$. Therefore $c \cdot \underset{e \geq e'}{\cap} \Ker dF^e_S =0$ for every $d \in S$ and $q'=p^{e'}$, namely, $c \in \tau(S,\Delta_S)$.
\end{proof}

Hara and Yoshida informed the author of the restriction theorem and subadditivity theorem for the ideal ``$\tau(R,\mathfrak{a})$'' (See \cite[Theorem 4.1 and Proposition 4.4]{HY}). Their results inspire him to obtain the restriction theorem and subadditivity theorem for ``$\tau(R,\Delta)$'' which correspond to Proposition \ref{DEL}.
\begin{prop}
\label{Hara}
\begin{enumerate}
\item \textup{(Restriction Theorem)} Let $(R,\mathfrak{m})$ be a complete normal Cohen-Macaulay $\mQ$-Gorenstein local ring of characteristic $p >0$, and assume that the order of the canonical class in the divisor class group $\mathrm{Cl}(R)$ is not divisible by $p$.
Let $\Delta$ be an effective $\mQ$-Cartier divisor on $\Spec R$, that is, $r\Delta = \Div_R(y)$ for some positive integer $r$ and nonzero element $y \in R$, and $x \in R$ be a nonzero divisor such that R/xR is normal and $y \notin xR$.  Then, denoting $S:=R/xR$, we have
$$\tau(S,\Delta|_{\Spec S}) \subseteq \tau(R,\Delta)\cdot S.$$
\item \textup{(Subadditivity Theorem)} Let $(R,\mathfrak{m})$ be a complete regular local ring of characteristic $p >0$, and $\Delta_1$ and $\Delta_2$ be any two effective $\mQ$-divisors on $\Spec R$. Then
$$\tau(R,\Delta_1+\Delta_2) \subseteq \tau(R,\Delta_1) \cdot \tau(R,\Delta_2).$$\end{enumerate}
\end{prop}
\begin{proof}
$(1)$ We identify $E_S$ with $(0:x)_{E_R}$.
\addtocounter{cln}{-2}
\begin{cln}
$$0^{*\Delta}_{E_R} \cap E_S \subseteq 0^{*\Delta|_{\Spec S}}_{E_S}.$$
\end{cln}
\begin{proof}[Proof of Claim $1$]
First we will look at the Frobenius actions on $E_R$ and $E_S$. Since $R$ is Cohen-Macaulay, we have the following commutative diagram with exact rows for infinitely many $q=p^e$ (See the proof of \cite[Theorem 4.9]{HW}). 
$$\xymatrix{0 \ar[r] & E_S \ar[r] \ar[d]^{F^e_S} & E_R \ar[d]^{x^{q-1}F^e_R} \\
0 \ar[r] & E_S \otimes_S {}^eS \ar[r] & E_R \otimes_R {}^eR
}$$
We fix sufficiently large such $q=p^e$.
If $\xi \in 0^{*\Delta}_{E_R} \cap E_S$, then for some nonzero element $c \in R$, we have $cF^e_R(\xi)=0$ in $E_R \otimes_R {}^eR((q-1)\Delta)$. 
We write $q-1=kr+i$ for integers $k$ and $i$ with $0 \leq i \leq r-1$.
Then there exists a nonzero element $c' \in R$ such that $c' \notin xR$ and $c'x^{q-1}y^{k}F^e_R(\xi)=0$ in $E_R \otimes_R {}^eR$. 
Hence, by the above diagram, $c'y^{k}F^e_S(\xi)=0$ in $E_S \otimes_S {}^eS$, so that $c'F^e_S(\xi)=0$ in $E_S \otimes_S {}^eS((q-1)\Delta|_{\Spec S})$ for infinitely many $q=p^e$. Since $c' \notin xR$, it implies that $\xi \in 0^{*\Delta|_{\Spec S}}_{E_S}$.
\end{proof}
Since $R$ is complete, by \cite[Lemma 3.3]{Ha3}, we have $0^{*\Delta}_{E_R} = (0:\tau(R,\Delta))_{E_R}$. 
Thus 
\begin{align*}
(0:x)_{0^{*\Delta}_{E_R}} &= (0:\tau(R,\Delta)+xR)_{E_R}
   = (0:\frac{\tau(R,\Delta)+xR}{xR})_{E_S} \\
  &=(0:\tau(R,\Delta)\cdot S)_{E_S}.
\end{align*}
In light of Claim $1$, $\tau(R,\Delta)\cdot S = \mathrm{Ann}_S(0:x)_{0^{*\Delta}_{E_R}} \supseteq \tau(S,\Delta|_{\Spec S})$.

$(2)$ First we consider the following claim.
\begin{cln}
Let $R=k[[x_1,\dots,x_n]]$ (resp. $S=k[[y_1,\dots,y_m]]$) be an $n$-dimensional (resp. $m$-dimensional) complete regular local ring over a field $k$ of characteristic $p>0$ and $\Delta_R$ (resp. $\Delta_S$) an effective $\mQ$-Weil divisor on $\Spec R$ (resp. $\Spec S$). Let $T = R \hat{\otimes}_k S = k[[x_1,\dots,x_n,y_1,\dots,y_m]]$, and we denote by $p_R:\Spec T \to \Spec R$ and $p_S:\Spec T \to \Spec S$ natural projections. Then
$$\tau(T,p_R^*\Delta_R+p_S^*\Delta_S) \subseteq (\tau(R,\Delta_R)\otimes_k \tau(S,\Delta_S))T.$$
\end{cln}
\begin{proof}[Proof of Claim $2$]
It suffices to show that
$$ 0^{*(p_R^*\Delta_R+p_S^*\Delta_S)}_{E_T} \supseteq 0^{*\Delta_R}_{E_R}\otimes_k E_S +E_R \otimes_k 0^{*\Delta_S}_{E_S},$$
but it is clear since $E_T = E_R \otimes_k E_S$.
\end{proof}
Let $\rho : T=R \hat{\otimes}_k R \to R$ be a diagonal map, and we denote by $p_1$ (resp. $p_2$) $:\Spec T \to \Spec R$ the first (resp. second) projection. 
Then the natural surjection $T \to T/\Ker \rho = R$ is a complete intersection, so it follows from the repeated application of the restriction theorem that $\tau(R,\Delta_1+\Delta_2) \subseteq \tau(T,p_1^*\Delta_1+p_2^*\Delta_2) \cdot R$. 
On the other hand, by Claim $2$, $\tau(T,p_1^*\Delta_1+p_2^*\Delta_2) \cdot R \subseteq \tau(R,\Delta_1) \cdot \tau(R,\Delta_2)$. 
Therefore the assertion follows.
\end{proof}

\begin{thm}\label{thm:test subset multi}
Let $(R, \mathfrak{m})$ be an F-finite normal local ring of characteristic $p > 0$ and $\Delta$ an effective $\mQ$-Weil divisor on $Y = \Spec R$ such that $K_Y + \Delta$ is $\mQ$-Cartier. Let $f:X \to Y=\Spec R$ be a proper birational morphism with $X$ normal. Then
$$\tau(R,\Delta) \subseteq H^0(X,\mathcal{O}_X(\lceil K_X -f^{*}(K_Y+\Delta) \rceil)).$$
\end{thm}
\begin{proof}
The essential idea of the proof is seen in \cite[Theorem 3.3]{HW} and \cite{Wa}. Our proof consists of six steps. 

(Step $1$) Take any nonzero element $c \in \tau(R,\Delta)$, and fix a nonzero element $d \in R(-\lceil \Delta \rceil)$. 
Then, for every $q'=p^{e'} > 0$, we have $c \cdot \underset{e \geq {e'}}{\cap} \Ker cdF^e = 0$, where $F^e:E_R \to E_R \otimes_R {}^eR((q-1)\Delta)$ is the $e$-times Frobenius map induced on $E_R$. Since $E_R$ is Artinian, there exists $q''=p^{e''}$ such that $\underset{e \geq {e'}}{\cap} \Ker cdF^e = \underset{{e''} \geq e \geq {e'}}{\cap} \Ker cdF^e$.
Let
$$ \varphi_e:\mathrm{Hom}_R(R((q-1)\Delta)^{1/q}, R) \to \mathrm{Hom}_R(R,R)=R $$
be an $R$-linear map induced by the $R$-linear map $R \xrightarrow{(cd)^{1/q}} R((q-1)\Delta)^{1/q}$, and set $\varphi = \underset{{e''} \geq e \geq {e'}}{\oplus}\varphi_e$. 
Since $cdF^e$ is the Matlis dual of $\varphi_e$, the condition that $c \cdot \underset{{e''} \geq e \geq {e'}}{\cap} \Ker cdF^e=0$ implies that $c \in \Img \varphi$. 
Hence, for every ${e''} \geq e \geq {e'}$, there exist $c_e \in R$ and an $R$-module homomorphism $ {\phi_e}':R((q-1)\Delta)^{1/q} \to R$ sending $(cd)^{1/q}$ to $c_e$ such that $\underset{{e''} \geq e \geq {e'}}{\sum}c_e=c$.

(Step $2$) We will prove that $\lfloor \Delta-\Div_R(c) \rfloor \leq 0$. Assume to the contrary that $\Delta$ has a component $\Delta_0$ such that the coefficient of $\Delta$ in $\Delta_0$ is at least $1+v_{\Delta_0}(c)$, where $v_{\Delta_0}$ is the valuation of $\Delta_0$. Since the coefficient of $(q-1)\Delta+\Div_R(d) \geq q\Delta$ in $\Delta_0$ is $q(1+v_{\Delta_0}(c))$ or more, the $R$-linear map $R \xrightarrow{(cd)^{1/q}} R((q-1)\Delta)^{1/q}$ factors through $R \hookrightarrow R((1+v_{\Delta_0}(c)){\Delta}_0)$. Hence, for every $e''\geq e \geq e'$, ${\phi_e}'$ induces an $R$-module homomorphism $R((1+v_{\Delta_0}(c)){\Delta}_0) \to R$ which sends $1$ to $c_e$. Thus there exists an $R$-linear map $R((1+v_{\Delta_0}(c)){\Delta}_0) \to R$ sending $1$ to $c$. This is a contradiction.

(Step $3$) Let $\phi_e= d^{1/q}{\phi_e}'$. 
The $R$-module homomorphism $\phi_e$ (resp. ${\phi_e}'$) induces an $R$-linear map $\psi_e$ (resp. ${\psi_e}'$) : $R((q-1)\Delta+\Div_R(c))^{1/q} \to R(\Div_R(c_e))$ which sends $1$ (resp. $d^{1/q}$) to $1$. 
We may assume without loss of generality that $X$ is Gorenstein (cf. the proof of \cite[Theorem 3.3]{HW}). 
Thanks to the adjunction formula, we may regard $\psi_e$ (resp. ${\psi_e}'$) in 
\begin{align*}
\mathrm{Hom}_R(R((q-1)\Delta+\Div_R(c))^{1/q}, R(\Div_R(c_e))) \\
&\hspace*{-2in}\cong R(\lceil(1-q)(K_Y+\Delta) +q\Div_R(c_e)-\Div_R(c) \rceil)^{1/q}
\end{align*}
as a rational section of the sheaf $\mathcal{O}_X((1-q)K_X)$, and consider the corresponding divisor on $X$
$$D_e = D_{\psi_e} = (\psi_e)_0 - (\psi_e)_{\infty} \quad (\text{resp. }{D_e}' = D_{{\psi_e}'} = ({\psi_e}')_0 - ({\psi_e}')_{\infty}),$$
where $(\psi_e)_0$ and  $(\psi_e)_{\infty}$ (resp. $({\psi_e}')_0$ and $({\psi_e}')_{\infty}$) are the divisors of zeros and poles of $\psi_e$ (resp. ${\psi_e}'$) as rational sections of $\mathcal{O}_X((1-q)K_X)$. 
Clearly, $D_e = {D_e}'+\Div_X(d)$. 
By definition, $D_e$ and ${D_e}'$ are linearly equivalent to $(1-q)K_X$, and $(\phi_e)_{\infty}$ and $({\phi_e}')_{\infty}$ are $f$-exceptional divisors. 
Hence $f_*D_e$ (resp. $f_*{D_e}'$) is linearly equivalent to $(1-q)K_Y$ and $f_*D_e$(resp. $f_*{D_e}'$) $\geq \lfloor (q-1)\Delta \rfloor -q\Div_R(c_e)+\Div_R(c)$. 
We denote $X' = X \setminus \text{Supp }(\psi_e)_{\infty}$. 
Then $\psi_e$ lies in
\begin{align*}
\mathrm{Hom}_{\mathcal{O}_{X'}}(\mathcal{O}_{X'}(\Div_{X'}(c))^{1/q},\mathcal{O}_{X'}(\Div_{X'}(c_e))) \\
&\hspace*{-2in}\cong  H^0(X',\mathcal{O}_{X'}((1-q)K_{X'}+q\Div_{X'}(c_e)-\Div_{X'}(c))).
\end{align*}

(Step $4$) We will show that the coefficient of $D_e$ in each irreducible component is $q-1$ or less. Assume to the contrary that there exists an irreducible component $D_{e,0}$ of $D_e$ whose coefficient is at least $q$. Let $v_{D_{e,0}}$ be the valuation of $D_{e,0}$ and $\alpha = qv_{D_{e,0}}(c_e)-v_{D_{e,0}}(c)+q$. 
Then $\psi_e$ lies in 
\begin{align*}
H^0(X', \mathcal{O}_{X'}(((1-q)K_{X'}+q\Div_{X'}(c_e)-\Div_{X'}(c))-\alpha D_{e,0})) \\
&\hspace*{-3in}\cong \mathrm{Hom}_{\mathcal{O}_{X'}}(\mathcal{O}_{X'}(\alpha D_{e,0}+\Div_{X'}(c))^{1/q},\mathcal{O}_{X'}(\Div_{X'}(c_e))).
\end{align*}
Therefore, we get the following commutative diagram. 
$$\xymatrix 
{\mathcal{O}_{X'} \ar@{^{(}->}[r] \ar@<-1ex>@{^{(}->}[dr] & \mathcal{O}_{X'}(q(v_{D_{e,0}}(c_e)+1)D_{e,0}+(\Div_{X'}(c)-v_{D_{e,0}}(c)D_{e,0}))^{1/q} \ar[d]^{\psi_e} \\
& \mathcal{O}_{X'}(\Div_{X'}(c_e)) \\
}$$
However the natural inclusion map
$$ \mathcal{O}_{X'} \hookrightarrow \mathcal{O}_{X'}(q(v_{D_{e,0}}(c_e)+1)D_{e,0}+(\Div_{X'}(c)-v_{D_{e,0}}(c)D_{e,0}))^{1/q} $$
factors through $\mathcal{O}_{X'}((v_{D_{e,0}}(c_e)+1)D_{e,0})$, and the above commutative diagram implies $(v_{D_{e,0}}(c_e)+1)D_{e,0} \leq \Div_{X'}(c_e)$. This is absurd. Hence every coefficient of $D_e$ in each irreducible component must be at most $q-1$.

(Step $5$) We denote by $\cup_{j=1}^s E_j$ the exceptional divisor of $f$ and by $f^{-1}_{*}\Delta'$ the strict transform of $\Delta':=\Delta-\Div_R(c)$ in $X$. Then we write
$$ K_X + f^{-1}_{*}\Delta' \underset{\text{$\mQ$-lin.}}{\sim} f^*(K_Y + \Delta') + \sum_{j=1}^{s}a_jE_j. $$
Let ${B_e}'=\frac{1}{q-1}{D_e}'-f^{-1}_{*}\Delta'$. Then ${B_e}'$ is $\mQ$-linearly equivalent to $-(K_X + f^{-1}_{*}\Delta')$, so that $f_*{{B_e}'}$ is $\mQ$-linearly equivalent to $-f_*(K_X + f^{-1}_{*}\Delta') = -(K_Y + \Delta')$.
Hence $f_*{B_e}'$ is $\mQ$-Cartier. Since ${B_e}' + \sum_{j=1}^s a_j E_j$ is $\mQ$-linearly equivalent to $-f^*(K_Y + \Delta')$, we know that $({B_e}' - f^*f_*{B_e}') + \sum_{j=1}^s a_j E_j$ is an $f$-exceptional divisor which is $\mQ$-linearly trivial relative to $f$. 
Hence
$$({B_e}' - f^*f_*{B_e}') + \sum_{j=1}^s a_j E_j = 0. $$

(Step $6$) Now, by (Step $3$),
\begin{align*}
f_*{D_e}' - (q-1)\Delta' &\geq (\lfloor (q-1)\Delta \rfloor -q\Div_R(c_e)+\Div_R(c))-(q-1)\Delta' \\
&\geq -{\Delta}''+q(\Div_R(c)-\Div_R(c_e))
\end{align*}
for some effective $\mQ$-Cartier divisor ${\Delta}''$ on $Y$ which is independent of $q$. 
This implies $f_*{B_e}' \geq -\frac{1}{q-1}{\Delta}''+\frac{q}{q-1}(\Div_R(c)-\Div_R(c_e))$, whence
$$f^*f_*B_e \geq - \frac{1}{q-1}f^*{\Delta}''+\frac{q}{q-1}(\Div_X(c)-\Div_X(c_e)).$$
On the other hand, we have seen in (Step $4$) that the coefficient of $D_e$ in $E_j$ is at most $q-1$. Since $D_e = {D_e}'+\Div_X(d)$ and we can assume that the coefficient of $\Div_X(d) - f^*{\Delta}''$ in $E_j$ is greater than zero, the coefficient of ${D_e}'+f^*\Delta''$ in $E_j$ is less than $q-1$.
Hence the coefficient of ${B_e}' - f^*f_*{B_e}'$ in $E_j$ is less than $1 - \frac{q}{q-1}(\Div_X(c)-\Div_X(c_e))$.
Since $\underset{{e''} \geq e \geq {e'}}{\sum}c_e=c$,
for every $j=1,\dots,s$, there exists $e'' \geq e \geq e'$ such that $v_{E_j}(c) \geq v_{E_j}(c_e)$, where $v_{E_j}$ is the valuation of $E_j$. Therefore via (Step $5$) we have $a_j > -1$.

It follows from the above result and (Step $2$) that
$$ \Div_X(c) + \lceil K_X-f^{*}(K_Y+\Delta)=\lceil K_X-f^{*}(K_Y+{\Delta}') \rceil \rceil \geq 0,$$
that is, $c \in H^0(X, \mathcal{O}_X(\lceil K_X -f^{*}(K_Y+\Delta) \rceil))$.
\end{proof}

\begin{rem}
By using the argument in (Step $1$) of the proof of Theorem \ref{thm:test subset multi}, we can also define the F-purity and F-regularity for ineffective divisors, and prove the similar results as those for effective divisors. However, the definition of F-purity and F-regularity for ineffective divisors is very complicated, therefore we only treat the effective case in this paper.
\end{rem}

\begin{expl}
\begin{enumerate}
\item When $R$ is a regular local ring and $\Delta$ is an effective $\mQ$-Weil divisor with simple normal crossing support, then $\tau(R,\Delta)$ $= R(-\lfloor \Delta \rfloor)$.
\item Let $R=k[[x_1,\dots,x_d]]$ be a $d$-dimensional complete regular local ring over a field $k$ of characteristic $p>0$ and $\Delta =\Div_R(x_1^{d+1}+\dots+x_d^{d+1})$. If the characteristic $p>d+1$, then $\tau(R, \frac{d}{d+1}\Delta) = (x_1,\dots,x_d)$. 
\end{enumerate}
\end{expl}

\section{Main Theorem}
To state the main result, we will explain the meaning of the phrase ``in characteristic $p \gg 0$.''

Let $R$ be a normal domain which is finitely generated over a field $k$ of characteristic zero and $\Delta$ an effective $\mQ$-Weil divisor on $\Spec R$ such that $K_Y + \Delta$ is $\mQ$-Cartier.
 Let $f:X \to \Spec R$ be a resolution of singularities such that $\mathrm{Exc}(f)+f^{-1}_*\Delta$ has simple normal crossing support. 
 Choosing a suitable finitely generated $\mZ$-subalgebra $A$ of $k$, there exists a finitely generated normal flat $A$-algebra $R_A$, an effective $\mQ$-Weil divisor $\Delta_A$ on $\Spec R_A$, a smooth $A$-scheme $X_A$ and a birational $A$-morphism $f_A:X_A \to \Spec R_A$ such that $K_{R_A} + \Delta_A$ is $\mQ$-Cartier, $\mathrm{Exc}(f_A)+{f_A}^{-1}_*\Delta_A$ has simple normal crossing support, and by tensoring $k$ over $A$ one gets back $R$, $\Delta$, $X$ and $f:X \to \Spec R$.
Given a closed point $s \in \Spec A$ with residue field $\kappa=\kappa(s)$, we denote the corresponding fibers over $s$ by $f_{\kappa}:X_\kappa \to \Spec R_\kappa$, $\Delta_\kappa$, etc. 
Then the pairs $(R_\kappa, \Delta_\kappa)$ over general closed points $s \in \Spec A$ inherit the properties possessed by the original one $(R, \Delta)$.

Now we fix a general closed point $s \in \Spec A$ with residue field $\kappa=\kappa(s)$ of sufficiently large characteristic $p \gg 0$. 
Then we refer to the fibers over $s \in \Spec A$ as ``reduction modulo $p \gg 0$,'' and use the phrase ``in characteristic $p \gg 0$'' when we look at general closed fibers which are reduced from characteristic zero to characteristic $p \gg 0$ as above.

The following lemma is essential to prove ``F-properties'' in characteristic $p \gg 0$.
\begin{lem}[\textup{\cite[Proposition 3.6, Corollary 3.8]{Ha2}}]\label{lem:injective}
Let $(R, \mathfrak{m})$ be a normal local ring of dimension $d \geq 2$, essentially of finite type over a perfect field $\kappa$ of characteristic $p > 0$. 
Let $f:X \to \Spec R$ be a resolution of singularities and $D$ an $f$-ample $\mQ$-Cartier divisor on $X$ with simple normal crossing support.
We denote the closed fiber of $f$ by $Z$. If $(R, \mathfrak{m})$ is the localization at any prime ideal of a finitely generated $\kappa$-algebra which is a reduction modulo $p \gg 0$ as well as $X, D$ and $f:X \to \Spec R$, then the $e$-times Frobenius map 
$$ F^e: H_{Z}^{d}(X, \mathcal{O} _X(-D)) \to H_Z^d (X, \mathcal{O} _X(-qD))$$
is injective for every $q=p^e$.
\end{lem}
Now we state our main result.
\begin{thm}\label{thm:mult=test}
Let $(R, \mathfrak{m})$ be a normal local ring essentially of finite type over a field of characteristic zero and $\Delta$ an effective $\mQ$-Weil divisor on $Y=\Spec R$ such that $K_Y+ \Delta$ is $\mQ$-Cartier. Then, in characteristic $p \gg 0$,
$$\tau (R, \Delta) = \mathcal{J}(Y, \Delta).$$
\end{thm}
\begin{proof}
Thanks to Theorem \ref{thm:test subset multi}, it suffices to prove that $\tau (R, \Delta) \supseteq \mathcal{J}(Y, \Delta)$ in characteristic $p \gg 0$.
Let $f:X \to Y=\Spec R$ be a resolution of singularities such that $\mathrm{Exc}(f)+f^{-1}_{*}\Delta$ has simple normal crossing support.

Take a nonzero element $c \in R(-\Delta_\mathrm{red})$ such that $R_c$ is regular, where $\Delta_\mathrm{red}$ is the reduced divisor whose support is equal to that of $\Delta$. Let $\Delta' = \Div_R(c)$. Then there is a rational number $0 \leq \epsilon \ll 1$ such that $ \lfloor f^*(K_Y + \Delta) \rfloor = \lfloor f^*(K_Y + \Delta+ \epsilon \Delta') \rfloor$. Take an $f$-ample $\mQ$-Cartier divisor $H$ on $X$ which is supported on the exceptional locus of $f$ such that $\lfloor f^*(K_Y+\Delta+\epsilon \Delta') - H\rfloor = \lfloor f^*(K_Y+\Delta)\rfloor$.
Set $D = H - f^*(K_Y+\Delta+\epsilon \Delta')$ and we may assume that $\mathrm{Exc}(f)+f^{-1}_{*}(\Delta+\epsilon \Delta')$ has simple normal crossing support again, replacing $f$ suitably.
By Lemma \ref{lem:injective}, in characteristic $p \gg0$, the $e$-times Frobenius map
$$ F^e:H_{Z}^{d}(X, \mathcal{O} _X(f^*(K_Y+\Delta))) \to H_Z^d (X, \mathcal{O} _X(-qD))$$
is injective for every $q=p^e$, where $Z$ is the closed fiber of $f$. 

On the other hand, let
$$\delta : H_{\mathfrak{m}}^d(R(K_Y)) \to H_{Z}^{d}(X, \mathcal{O} _X(f^*(K_Y+\Delta)))$$ 
be the Matlis dual of the natural inclusion map 
$$ \mathcal{J}(Y, \Delta) = H^0(X, \mathcal{O}_X(\lceil K_X-f^*(K_Y+\Delta)\rceil)) \hookrightarrow R,$$
and 
$$\delta_e:H_{\mathfrak{m}}^d(R(q(K_Y+\Delta + \epsilon \Delta'))) \to H_Z^d (X, \mathcal{O} _X(-qH+qf^{*}(K_Y+\Delta+\epsilon \Delta')))$$
the natural map induced by an edge map of the Leray spectral sequence
$$H^j_{\mathfrak{m}}(H^i(X,\mathcal{O}_X(qf^{*}(K_Y+\Delta+\epsilon \Delta')))) \Rightarrow H_Z^{i+j} (X, \mathcal{O} _X(qf^*(K_Y+\Delta+\epsilon \Delta'))).$$
Then by the Matlis duality,
\begin{align*}
 \ker(\delta)
  & = \mathrm{Hom}_R \left( \frac{R}{\mathcal{J}(Y, \Delta)}\,,\,E_R \right)
  = \mathrm{Ann}_{H_{\mathfrak{m}}^d(R(K_Y))}\mathcal{J}(Y, \Delta), \\
 \ker(\delta_e) & = \mathrm{Hom}_R \left( \frac{R(\lceil K_Y-q(K_Y+\Delta+\epsilon \Delta')\rceil)}{H^0(X, \mathcal{O}_X(\lceil K_X + qD \rceil))} \,,\,E_R \right) \\
 & =\mathrm{Ann}_{H_{\mathfrak{m}}^d(R(q(K_Y+\Delta+\epsilon \Delta')))}H^0(X, \mathcal{O}_X(\lceil K_X + qD \rceil)).
\end{align*}

We obtain the following commutative diagram with exact rows for every $q=p^e$.
$$\xymatrix{
0 \to \ker(\delta) \ar[d] \ar[r] & H_{\mathfrak{m}}^d(R(K_Y)) \ar[d]^{F^e} \ar[r]^{\delta \qquad} & H_Z^d (X, \mathcal{O} _X(-D)) \ar[d]^{F^e} \to 0 \\
0 \to \ker(\delta_e) \ar[r] & H_{\mathfrak{m}}^d(R(q(K_Y+\Delta + \epsilon \Delta'))) \ar[r]^{\quad \delta_e} & H_Z^d (X, \mathcal{O} _X(-qD)) \to 0 \\
}$$
Take any element $\xi \in H_{\mathfrak{m}}^d(R(K_Y)) \setminus \ker(\delta)$. By the above diagram, ${\xi}^q \notin \ker(\delta_e)$. By Lemma \ref{lem:test-element}, $c^n$ is a $\Delta$-test element for some positive integer $n$, and for sufficiently large $q$, 
\begin{align*}
H^0(X, \mathcal{O}_X(\lceil K_X + q(H-f^*(K_Y+\Delta+\epsilon \Delta')) \rceil)) \\
&\hspace*{-2in}\subseteq c^{n+1} H^0(X, \mathcal{O}_X(\lceil K_X+q(H-f^*(K_Y+\Delta)) \rceil)).
\end{align*}
Hence
$$c^{n+1} \xi^q \notin \mathrm{Ann}_{H_{\mathfrak{m}}^d(R(q(K_Y+\Delta)))}H^0(X, \mathcal{O}_X(\lceil K_X + q(H-f^*(K_Y+\Delta)) \rceil)).$$
If $\xi \in 0_{H_{\mathfrak{m}}^d(R(K_Y))}^{*\Delta}$, then by the proof of Lemma \ref{lem:test-element}, we have $c^n \xi^q = 0$ in $H_{\mathfrak{m}}^d(R(\lceil q K_Y+ q\Delta \rceil))$. 
Therefore $c^{n+1} \xi^q = 0$ in $H_{\mathfrak{m}}^d(R(q (K_Y+\Delta)))$, and this is a contradiction. 
It follows that 
$$0_{H_{\mathfrak{m}}^d(R(K_Y))}^{*\Delta} \subseteq \ker(\delta) = \mathrm{Ann}_{H_{\mathfrak{m}}^d(R(K_Y))}\mathcal{J}(Y, \Delta),$$
and by Matlis duality (see \cite[Lemma 3.3]{Ha3}), $\tau (R, \Delta) \supseteq \mathcal{J}(Y, \Delta)$.
\end{proof}
\begin{rem}
When $\Delta=0$, Theorem \ref{thm:mult=test} coincides with the results of Hara \cite{Ha3} and Smith \cite{Sm2}. However, since $R$ is not necessarily $\mQ$-Gorenstein in our situation, we cannot use their strategy which is to reduce the case where $R$ is quasi-Gorenstein by passing to an index one cover.
\end{rem}

As a direct consequence of the main theorem, we get the equivalence of klt pairs and strongly F-regular pairs.
\begin{cor}[\textup{\cite[Conjecture 5.1.1]{HW}}]\label{cor:F-reg=klt}
Let $(R, \mathfrak{m})$ be a normal local ring essentially of finite type over a field of characteristic zero and $\Delta$ an effective $\mQ$-Weil divisor on $Y=\Spec R$ such that $K_Y+ \Delta$ is $\mQ$-Cartier. Then, $(Y, \Delta)$ is klt if and only if $(R, \Delta)$ is of strongly F-regular type.
\end{cor}

Hara and K.-i.~Watanabe \cite[Problem 5.1.2]{HW} conjectured that $(Y, \Delta)$ is lc if and only if $(R, \Delta)$ is of dense F-pure type.
The following result about log canonical thresholds is a piece of evidence for their conjecture. 
See \cite{Ko} for the basic properties of log canonical thresholds.
\begin{cor}[\textup{\cite[Conjecture 5.2.1]{HW}}]
Let $Y$ be a normal variety in characteristic zero with only klt singularity at a point $y \in Y$ and $\Delta$ be an effective $\mQ$-Cartier divisor on $Y$. 
We denote by $C_y(Y, \Delta)$ the log canonical threshold of $\Delta$ at $y \in Y$, that is, 
\begin{align*}
C_y(Y, \Delta) &= \sup \{t \in \mR_{>0} \mid (Y, t\Delta) \textup{ is lc at }y \in Y \}\\
               &= \sup \{t \in \mR_{>0} \mid (Y, t\Delta) \textup{ is klt at }y \in Y \}\text{.}
\end{align*}
Then, 
\begin{align*}
C_y(Y, \Delta) &= \sup \{t \in \mQ_{>0} \mid (\mathcal{O}_{Y,y}, t\Delta) \textup{ is of dense F-pure type} \} \\ 
               &= \sup \{t \in \mQ_{>0} \mid (\mathcal{O}_{Y,y}, t\Delta) \textup{ is of strongly F-regular type} \} \text{.}
\end{align*}
\end{cor}
\begin{proof}
By \cite[Theorem 3.7]{HW}, the pair of dense F-pure type is lc. Hence the assertion is clear.
\end{proof}

\end{document}